\documentclass[11pt,a4paper]{article}
\usepackage{theorem,enumerate}
\usepackage{amsmath,latexsym,amssymb,amsfonts}
\usepackage{eucal}
\usepackage{color}
\usepackage{comment}
\newcounter{Scounter}
\theorembodyfont{\normalfont\slshape}
\newtheorem{thm}{Theorem}[section]

\newtheorem{Thm}{Theorem}

\newtheorem{prop}[thm]{Proposition}
\newtheorem{lem}[thm]{Lemma}

\numberwithin{equation}{section}
\newtheorem{remark}{Remark}

\newcommand{\proof}{\medbreak\noindent\textit{Proof.}\quad}

\newcommand{\qed}{{$\quad\square$\vs{3.6}}}

\newcommand{\vs}[1]{\vspace*{#1 mm}}

\bfseries\normalfont

\title{A continuous generalization of domination-like invariants}
\author{
Michitaka Furuya\footnote{\texttt{e-mail:michitaka.furuya@gmail.com}} \vs{5}\\
\textsl{College of Liberal Arts and Science,}\\
\textsl{Kitasato University,}\\
\textsl{1-15-1 Kitasato, Minami-ku, Sagamihara, Kanagawa 252-0373, Japan}
}

\date{}

\begin{document}

\maketitle

\begin{abstract}
In this paper, we define a new domination-like invariant of graphs.
Let $\mathbb{R}^{+}$ be the set of non-negative numbers.
Let $c\in \mathbb{R}^{+}-\{0\}$ be a number, and let $G$ be a graph.
A function $f:V(G)\rightarrow \mathbb{R}^{+}$ is a $c$-self-dominating function of $G$ if for every $u\in V(G)$, $f(u)\geq c$ or $\max\{f(v):v\in N_{G}(u)\}\geq 1$.
The $c$-self-domination number $\gamma ^{c}(G)$ of $G$ is defined as $\gamma ^{c}(G):=\min\{\sum _{u\in V(G)}f(u):f$ is a $c$-self-dominating function of $G\}$.
Then $\gamma ^{1}(G)$, $\gamma ^{\infty }(G)$ and $\gamma ^{\frac{1}{2}}(G)$ are equal to the domination number, the total domination number and the half of the Roman domination number of $G$, respectively.
Our main aim is to continuously fill in the gaps among such three invariants.
In this paper, we give a sharp upper bound of the $c$-self-domination number for all $c\geq \frac{1}{2}$.
\end{abstract}



%%%%%%%%%%%%%%%%%%%%%%%%%%%%%%%%%%%%%%%%%%%%%%%%%%%%%%%%%%%%%%%%%%%%%%%%%%%%%%%%%%%%%%%%%%%%%%%%%%%%%%%%%%%%%%%%%%%%%%%%
%%%%%%%%%%%%%%%%%%%%%%%%%%%%%%%%%%%%%%%%%%%%%%%%%%%%%%%%%%%%%%%%%%%%%%%%%%%%%%%%%%%%%%%%%%%%%%%%%%%%%%%%%%%%%%%%%%%%%%%%
%%%%%%%%%%%%%%%%%%%%%%%%%%%%%%%%%%%%%%%%%%%%%%%%%%%%%%%%%%%%%%%%%%%%%%%%%%%%%%%%%%%%%%%%%%%%%%%%%%%%%%%%%%%%%%%%%%%%%%%%
\section{Introduction}\label{sec1}
%%%%%%%%%%%%%%%%%%%%%%%%%%%%%%%%%%%%%%%%%%%%%%%%%%%%%%%%%%%%%%%%%%%%%%%%%%%%%%%%%%%%%%%%%%%%%%%%%%%%%%%%%%%%%%%%%%%%%%%%
%%%%%%%%%%%%%%%%%%%%%%%%%%%%%%%%%%%%%%%%%%%%%%%%%%%%%%%%%%%%%%%%%%%%%%%%%%%%%%%%%%%%%%%%%%%%%%%%%%%%%%%%%%%%%%%%%%%%%%%%
%%%%%%%%%%%%%%%%%%%%%%%%%%%%%%%%%%%%%%%%%%%%%%%%%%%%%%%%%%%%%%%%%%%%%%%%%%%%%%%%%%%%%%%%%%%%%%%%%%%%%%%%%%%%%%%%%%%%%%%%

%%%%%%%%%%%%%%%%%%%%%%%%%%%%%%%%%%%%%%%%%%%%%%%%%%%%%%%%%%%%%%%%%%%%%%%%%%%%%%%%%%%%%%%%%%%%%%%%%%%%%%%%%%%%%%%%%%%%%%%%
%%%%%%%%%%%%%%%%%%%%%%%%%%%%%%%%%%%%%%%%%%%%%%%%%%%%%%%%%%%%%%%%%%%%%%%%%%%%%%%%%%%%%%%%%%%%%%%%%%%%%%%%%%%%%%%%%%%%%%%%
\subsection{Definitions and notations}\label{sec1.1}
%%%%%%%%%%%%%%%%%%%%%%%%%%%%%%%%%%%%%%%%%%%%%%%%%%%%%%%%%%%%%%%%%%%%%%%%%%%%%%%%%%%%%%%%%%%%%%%%%%%%%%%%%%%%%%%%%%%%%%%%
%%%%%%%%%%%%%%%%%%%%%%%%%%%%%%%%%%%%%%%%%%%%%%%%%%%%%%%%%%%%%%%%%%%%%%%%%%%%%%%%%%%%%%%%%%%%%%%%%%%%%%%%%%%%%%%%%%%%%%%%

All graphs considered in this paper are finite, simple, and undirected.
Let $G$ be a graph.
We let $V(G)$ and $E(G)$ denote the {\it vertex set} and the {\it edge set} of $G$, respectively.
For a vertex $u\in V(G)$, we let $N_{G}(u)$ and $d_{G}(u)$ denote the {\it neighborhood} and the {\it degree} of $u$, respectively; thus $N_{G}(u)=\{v\in V(G): uv\in E(G)\}$ and $d_{G}(u)=|N_{G}(u)|$.
For a subset $U$ of $V(G)$, we let $G[U]$ denote the subgraph of $G$ induced by $U$.
We let $P_{n}$ denote the {\it path} of order $n$.
For terms and symbols not defined in this paper, we refer the reader to \cite{D}.

We let $\mathbb{R}^{+}$ denote the set of non-negative numbers.
Here we regard $\infty $ as a non-negative number (i.e., $\infty \in \mathbb{R}^{+}$).
For a graph $G$ and a function $f:V (G)\rightarrow \mathbb{R}^{+}$, the {\it weight} $w(f)$ of $f$ is defined by $w(f)=\sum _{u\in V(G)}f(u)$.

Let $G$ be a graph.
A set $S\subseteq V(G)$ is a {\it dominating set} of $G$ if each vertex in $V(G)-S$ is adjacent to a vertex in $S$.
The minimum cardinality of a dominating set of $G$, denoted by $\gamma (G)$, is called the {\it domination number} of $G$.
The domination number is one of the most important invariants in Graph Theory, and it can be widely applied to real problems, for example, school bus routing problem, social network theory and location of radio stations (see \cite{HHS1,HHS2}).
To meet various additional requirements for above problems, many domination-like invariants were defined and studied.

A set $S\subseteq V(G)$ is a {\it total dominating set} of $G$ if each vertex of $G$ is adjacent to a vertex in $S$.
Note that if $G$ has no isolated vertices, then there exists a total dominating set of $G$.
For a graph $G$ without isolated vertices, the minimum cardinality of a total dominating set of $G$, denoted by $\gamma _{t}(G)$, is called the {\it total domination number} of $G$.
The total domination number is typically defined for only graphs without isolated vertices.
However, in this paper, we define $\gamma _{t}(G)$ as $\gamma _{t}(G)=\infty $ if $G$ has an isolated vertex for convenience.
The concept of total domination was introduced in \cite{CDH}, and has been actively studied (see a book~\cite{HY}).

A function $f:V (G)\rightarrow \{0,1,2\}$ is a {\it Roman dominating function} of $G$ if each vertex $u\in V(G)$ with $f(u)=0$ is adjacent to a vertex $v\in V(G)$ with $f(v)=2$.
The minimum weight of a Roman dominating function of $G$, denoted by $\gamma _{R}(G)$, is called the {\it Roman domination number} of $G$.
The Roman domination number was introduced by Stewart~\cite{S}, and was studied by Cockayne et al.~\cite{CDHH} in earnest.
Roman domination derives from the strategy to defend the Roman Empire against the enemies.
Recently, various properties on the Roman domination number have been explored in, for example, \cite{FYJ,LKLP,LC}.

From a mathematical point of view, Roman domination concept seems to be more artificial than original domination and total domination.
However, by the following reasons, we can interpret Roman domination as a natural extension of domination and total domination.

We define a new domination-like invariant.
Let $c\in \mathbb{R}^{+}-\{0\}$ be a number, and let $G$ be a graph.
A function $f:V(G)\rightarrow \mathbb{R}^{+}$ is a {\it $c$-self-dominating function} (or {\it $c$-SDF}) of $G$ if for each $u\in V(G)$, $f(u)\geq c$ or $\max\{f(v):v\in N_{G}(u)\}\geq 1$.

\begin{remark}%%%%%%%%%%%%%%%%%%%%%%%%%%%%%%%%%%%%%%%%%%%%%%%%%%%%%%%%%%%%%%%%%%%%%%%%%%%%%%%%%%%%%%%%%%%%%%%%%%%%%%%%%%
\label{remark1}
We choose a $c$-SDF $f$ of $G$ so that $w(f)$ is as small as possible.
Suppose that there exists a vertex $x\in V(G)$ with $f(x)\notin \{0,1,c\}$.
Now we construct a function $f':V(G)\rightarrow \mathbb{R}^{+}$ as follows:
For $u\in V(G)-\{x\}$, let $f'(u)=f(u)$.
If $0<f(x)<\min\{1,c\}$, let $f'(x)=0$; if $\min\{1,c\}<f(x)<\max\{1,c\}$, let $f'(x)=\min\{1,c\}$; if $f(x)>\max\{1,c\}$, let $f'(x)=\max\{1,c\}$ (for example, if $c=\infty $, then $f'(x)=0$ or $f'(x)=1$ according as $0<f(x)<1$ or $1<f(x)<\infty $).
Then we can easily verify that $f'$ is a $c$-SDF of $G$ with $w(f')<w(f)$, which contradicts the choice of $f$.
Thus $\{f(u):u\in V(G)\}\subseteq \{0,1,c\}$.
In particular, the minimum weight of $c$-SDF of $G$ is well-defined.
\end{remark}
%%%%%%%%%%%%%%%%%%%%%%%%%%%%%%%%%%%%%%%%%%%%%%%%%%%%%%%%%%%%%%%%%%%%%%%%%%%%%%%%%%%%%%%%%%%%%%%%%%%%%%%%%%%%%%%%%%%%%%%%

The minimum weight of a $c$-SDF of $G$, denoted by $\gamma ^{c}(G)$, is called the {\it $c$-self-domination number} of $G$.
A $c$-SDF $f$ of $G$ with $w(f)=\gamma ^{c}(G)$ is called a {\it $\gamma ^{c}$-function} of $G$.
Considering Remark~\ref{remark1},
\begin{align}
\mbox{$\{f(u):u\in V(G)\}\subseteq \{0,1,c\}$ for a $\gamma ^{c}$-function $f$ of a graph $G$.}\label{01c}
\end{align}

We show that $c$-self-domination is a common generalization of domination, total domination and Roman domination.

\begin{prop}%%%%%%%%%%%%%%%%%%%%%%%%%%%%%%%%%%%%%%%%%%%%%%%%%%%%%%%%%%%%%%%%%%%%%%%%%%%%%%%%%%%%%%%%%%%%%%%%%%%%%%%%%%%%
\label{prop1.1}
Let $G$ be a graph.
Then the following hold.
\begin{enumerate}
\item[{\upshape(i)}]
$\gamma ^{1}(G)=\gamma (G)$,
\item[{\upshape(ii)}]
$\gamma ^{\infty }(G)=\gamma _{t}(G)$, and
\item[{\upshape(iii)}]
$\gamma ^{\frac{1}{2}}(G)=\frac{1}{2}\gamma _{R}(G)$.
\end{enumerate}
\end{prop}
%%%%%%%%%%%%%%%%%%%%%%%%%%%%%%%%%%%%%%%%%%%%%%%%%%%%%%%%%%%%%%%%%%%%%%%%%%%%%%%%%%%%%%%%%%%%%%%%%%%%%%%%%%%%%%%%%%%%%%%%
\proof
\begin{enumerate}
\item[{\upshape(i)}]
For a dominating set $S$ of $G$ with $|S|=\gamma (G)$, the function $f_{1}:V(G)\rightarrow \mathbb{R}^{+}$ with
$$
f_{1}(u)=
\begin{cases}
1, & u\in S;\\
0, & u\notin S
\end{cases}
$$
is a $1$-SDF of $G$ with $w(f_{1})=|S|$, and hence $\gamma ^{1}(G)\leq w(f_{1})=|S|=\gamma (G)$.

Let $f$ be a $\gamma ^{1}$-function of $G$.
Then by (\ref{01c}), $\{f(u):u\in V(G)\}\subseteq \{0,1\}$.
Hence the set $S_{1}:=\{u\in V(G):f(u)=1\}$ is a dominating set of $G$ with $|S_{1}|=w(f)$.
Thus $\gamma (G)\leq |S_{1}|=w(f)=\gamma ^{1}(G)$.

Consequently, $\gamma ^{1}(G)=\gamma (G)$.

\item[{\upshape(ii)}]
If $G$ has an isolated vertex, then it is clear that $\gamma ^{\infty }(G)=\infty =\gamma _{t}(G)$.
Thus we may assume that $G$ has no isolated vertices.
For a total dominating set $S$ of $G$ with $|S|=\gamma _{t}(G)$, the function $f_{2}:V(G)\rightarrow \mathbb{R}^{+}$ with
$$
f_{2}(u)=
\begin{cases}
1, & u\in S;\\
0, & u\notin S
\end{cases}
$$
is an $\infty $-SDF of $G$ with $w(f_{2})=|S|$, and hence $\gamma ^{\infty }(G)\leq w(f_{2})=|S|=\gamma _{t}(G)$.

Let $f$ be a $\gamma ^{\infty }$-function of $G$.
Since the function assigning $1$ to all vertices of $G$ is an $\infty $-SDF of $G$, we have $\gamma ^{\infty }(G)<\infty $, and hence $f$ does not use $\infty $.
This together with (\ref{01c}) implies that $\{f(u):u\in V(G)\}\subseteq \{0,1\}$.
Hence the set $S_{2}:=\{u\in V(G):f(u)=1\}$ is a total dominating set of $G$ with $|S_{2}|=w(f)$.
Thus $\gamma _{t}(G)\leq |S_{2}|=w(f)=\gamma ^{\infty }(G)$.

Consequently, $\gamma ^{\infty }(G)=\gamma _{t}(G)$.

\item[{\upshape(iii)}]
For a Roman dominating function $f$ of $G$ with $w(f)=\gamma _{R}(G)$, the function $f_{3}:V(G)\rightarrow \mathbb{R}^{+}$ with $f_{3}(u)=\frac{1}{2}f(u)~(u\in V(G))$ is a $\frac{1}{2}$-SDF of $G$ with $w(f_{3})=\frac{1}{2}w(f)$, and hence $\gamma ^{\frac{1}{2}}(G)\leq w(f_{3})=\frac{1}{2}w(f)=\frac{1}{2}\gamma _{R}(G)$.

Let $f'$ be a $\gamma ^{\frac{1}{2}}$-function of $G$.
Then by (\ref{01c}), $\{f'(u):u\in V(G)\}\subseteq \{0,1,\frac{1}{2}\}$.
Hence the function $f'_{3}:V(G)\rightarrow \{0,1,2\}$ with $f'_{3}(u)=2f'(u)~(u\in V(G))$ is a Roman dominating function of $G$ with $w(f'_{3})=2w(f')$.
Thus $\gamma _{R}(G)\leq w(f'_{3})=2w(f')=2\gamma ^{\frac{1}{2}}(G)$.

Consequently, $\gamma ^{\frac{1}{2}}(G)=\frac{1}{2}\gamma _{R}(G)$.
\qed
\end{enumerate}

%%%%%%%%%%%%%%%%%%%%%%%%%%%%%%%%%%%%%%%%%%%%%%%%%%%%%%%%%%%%%%%%%%%%%%%%%%%%%%%%%%%%%%%%%%%%%%%%%%%%%%%%%%%%%%%%%%%%%%%%
%%%%%%%%%%%%%%%%%%%%%%%%%%%%%%%%%%%%%%%%%%%%%%%%%%%%%%%%%%%%%%%%%%%%%%%%%%%%%%%%%%%%%%%%%%%%%%%%%%%%%%%%%%%%%%%%%%%%%%%%
\subsection{Main results}\label{sec1.2}
%%%%%%%%%%%%%%%%%%%%%%%%%%%%%%%%%%%%%%%%%%%%%%%%%%%%%%%%%%%%%%%%%%%%%%%%%%%%%%%%%%%%%%%%%%%%%%%%%%%%%%%%%%%%%%%%%%%%%%%%
%%%%%%%%%%%%%%%%%%%%%%%%%%%%%%%%%%%%%%%%%%%%%%%%%%%%%%%%%%%%%%%%%%%%%%%%%%%%%%%%%%%%%%%%%%%%%%%%%%%%%%%%%%%%%%%%%%%%%%%%

By Proposition~\ref{prop1.1}, $c$-self-domination can continuously fill in the gaps among three invariants: domination, total domination and Roman domination.
On the other hand, some results concerning such invariants have been proved via different techniques.
Thus the study of $c$-self-domination for $c\geq \frac{1}{2}$ may give essential boundaries of them.
As the initial research for the goal, we focus on the following known upper bounds.

\begin{Thm}[Ore~\cite{O}]%%%%%%%%%%%%%%%%%%%%%%%%%%%%%%%%%%%%%%%%%%%%%%%%%%%%%%%%%%%%%%%%%%%%%%%%%%%%%%%%%%%%%%%%%%%%%%%
\label{ThmA}
Let $G$ be a connected graph of order $n\geq 2$.
Then $\gamma (G)\leq \frac{1}{2}n$.
\end{Thm}
%%%%%%%%%%%%%%%%%%%%%%%%%%%%%%%%%%%%%%%%%%%%%%%%%%%%%%%%%%%%%%%%%%%%%%%%%%%%%%%%%%%%%%%%%%%%%%%%%%%%%%%%%%%%%%%%%%%%%%%%

\begin{Thm}[Cockayne et al.~\cite{CDH}]%%%%%%%%%%%%%%%%%%%%%%%%%%%%%%%%%%%%%%%%%%%%%%%%%%%%%%%%%%%%%%%%%%%%%%%%%%%%%%%%%
\label{ThmB}
Let $G$ be a connected graph of order $n\geq 3$.
Then $\gamma _{t}(G)\leq \frac{2}{3}n$.
\end{Thm}
%%%%%%%%%%%%%%%%%%%%%%%%%%%%%%%%%%%%%%%%%%%%%%%%%%%%%%%%%%%%%%%%%%%%%%%%%%%%%%%%%%%%%%%%%%%%%%%%%%%%%%%%%%%%%%%%%%%%%%%%

\begin{Thm}[Chambers et al.~\cite{CKPW}]%%%%%%%%%%%%%%%%%%%%%%%%%%%%%%%%%%%%%%%%%%%%%%%%%%%%%%%%%%%%%%%%%%%%%%%%%%%%%%%%
\label{ThmC}
Let $G$ be a connected graph of order $n\geq 3$.
Then $\gamma _{R}(G)\leq \frac{4}{5}n$.
\end{Thm}
%%%%%%%%%%%%%%%%%%%%%%%%%%%%%%%%%%%%%%%%%%%%%%%%%%%%%%%%%%%%%%%%%%%%%%%%%%%%%%%%%%%%%%%%%%%%%%%%%%%%%%%%%%%%%%%%%%%%%%%%

In this paper, we generalize Theorems~\ref{ThmA}--\ref{ThmC} as follows.

\begin{thm}%%%%%%%%%%%%%%%%%%%%%%%%%%%%%%%%%%%%%%%%%%%%%%%%%%%%%%%%%%%%%%%%%%%%%%%%%%%%%%%%%%%%%%%%%%%%%%%%%%%%%%%%%%%%%
\label{mainthm}
Let $c\in \mathbb{R}^{+}$ be a number with $c\geq \frac{1}{2}$.
Let $G$ be a connected graph of order $n\geq 3$.
Then
$$
\gamma ^{c}(G)\leq 
\begin{cases}
\frac{m+1}{2m+3}n, & \frac{m}{m+1}\leq c<\frac{2m+1}{2m+3},~m\in \mathbb{N};\\
\frac{cm+2c+1}{2m+5}n, & \frac{2m+1}{2m+3}\leq c<\frac{m+1}{m+2},~m\in \mathbb{N};\\
\frac{1}{2}n, & c=1;\\
\frac{m+2}{2m+3}n, & \frac{m+2}{m+1}<c\leq \frac{(2m+1)(m+2)}{(2m+3)m},~m\in \mathbb{N};\\
\frac{cm}{2m+1}n, & \frac{(2m+1)(m+2)}{(2m+3)m}<c\leq \frac{m+1}{m},~m\in \mathbb{N};\\
\frac{2}{3}n, & c>2.
\end{cases}
$$
\end{thm}
%%%%%%%%%%%%%%%%%%%%%%%%%%%%%%%%%%%%%%%%%%%%%%%%%%%%%%%%%%%%%%%%%%%%%%%%%%%%%%%%%%%%%%%%%%%%%%%%%%%%%%%%%%%%%%%%%%%%%%%%

\begin{remark}%%%%%%%%%%%%%%%%%%%%%%%%%%%%%%%%%%%%%%%%%%%%%%%%%%%%%%%%%%%%%%%%%%%%%%%%%%%%%%%%%%%%%%%%%%%%%%%%%%%%%%%%%%
\begin{enumerate}
\item[{\upshape(i)}]
Since the function $h(m):=\frac{m}{m+1}$ is a monotonically increasing function and $\lim_{m\rightarrow \infty }h(m)=1$, the interval $[\frac{1}{2},1)$ can be partitioned by intervals $[\frac{m}{m+1},\frac{m+1}{m+2})$ $(m\in \mathbb{N})$.
In particular, for a number $c~(\frac{1}{2}\leq c<1)$, there is only one positive integer $m$ such that $\frac{m}{m+1}\leq c<\frac{m+1}{m+2}$.
Furthermore, since $\frac{m}{m+1}<\frac{2m+1}{2m+3}<\frac{m+1}{m+2}$, the interval $[\frac{m}{m+1},\frac{m+1}{m+2})$ can be partitioned by two intervals $[\frac{m}{m+1},\frac{2m+1}{2m+3})$ and $[\frac{2m+1}{2m+3},\frac{m+1}{m+2})$.

\item[{\upshape(ii)}]
Since the function $h'(m):=\frac{m+1}{m}$ is a monotonically decreasing function and $\lim_{m\rightarrow \infty }h'(m)=1$, the interval $(1,2]$ can be partitioned by intervals $(\frac{m+2}{m+1},\frac{m+1}{m}]$ $(m\in \mathbb{N})$.
In particular, for a number $c~(1<c\leq 2)$, there is only one positive integer $m$ such that $\frac{m+2}{m+1}<c\leq \frac{m+1}{m}$.
Furthermore, since $\frac{m+2}{m+1}<\frac{(2m+1)(m+2)}{(2m+3)m}<\frac{m+1}{m}$, the interval $(\frac{m+2}{m+1},\frac{m+1}{m}]$ can be partitioned by two intervals $(\frac{m+2}{m+1},\frac{(2m+1)(m+2)}{(2m+3)m}]$ and $(\frac{(2m+1)(m+2)}{(2m+3)m},\frac{m+1}{m}]$.
\end{enumerate}

\noindent
By (i) and (ii), for each $c~(c\geq \frac{1}{2})$, Theorem~\ref{mainthm} gives exactly one upper bound on $\gamma ^{c}$.
\end{remark}
%%%%%%%%%%%%%%%%%%%%%%%%%%%%%%%%%%%%%%%%%%%%%%%%%%%%%%%%%%%%%%%%%%%%%%%%%%%%%%%%%%%%%%%%%%%%%%%%%%%%%%%%%%%%%%%%%%%%%%%%

It follows from Proposition~\ref{prop1.1}(i) and Theorem~\ref{ThmA} that $\gamma ^{1}(G)\leq \frac{1}{2}n$ for every connected graph $G$ of order $n\geq 2$.
(Indeed, since a maximal independent set $S$ of $G$ and the set $S':=V(G)-S$ are dominating sets of $G$, we obtain the upper bound.)
Thus it suffices to focus on Theorem~\ref{mainthm} for the case where $c\neq 1$.
We divide the proof of Theorem~\ref{mainthm} into three cases.
We consider the case where $\frac{1}{2}\leq c<1$ in Section~\ref{sec3}, the case where $1<c\leq 2$ in Section~\ref{sec4}, the case where $c>2$ in Section~\ref{sec5}.
In Section~\ref{sec6}, we discuss the sharpness of Theorem~\ref{mainthm}.

%%%%%%%%%%%%%%%%%%%%%%%%%%%%%%%%%%%%%%%%%%%%%%%%%%%%%%%%%%%%%%%%%%%%%%%%%%%%%%%%%%%%%%%%%%%%%%%%%%%%%%%%%%%%%%%%%%%%%%%%
%%%%%%%%%%%%%%%%%%%%%%%%%%%%%%%%%%%%%%%%%%%%%%%%%%%%%%%%%%%%%%%%%%%%%%%%%%%%%%%%%%%%%%%%%%%%%%%%%%%%%%%%%%%%%%%%%%%%%%%%
%%%%%%%%%%%%%%%%%%%%%%%%%%%%%%%%%%%%%%%%%%%%%%%%%%%%%%%%%%%%%%%%%%%%%%%%%%%%%%%%%%%%%%%%%%%%%%%%%%%%%%%%%%%%%%%%%%%%%%%%
\section{Trees without good edges}\label{sec2}
%%%%%%%%%%%%%%%%%%%%%%%%%%%%%%%%%%%%%%%%%%%%%%%%%%%%%%%%%%%%%%%%%%%%%%%%%%%%%%%%%%%%%%%%%%%%%%%%%%%%%%%%%%%%%%%%%%%%%%%%
%%%%%%%%%%%%%%%%%%%%%%%%%%%%%%%%%%%%%%%%%%%%%%%%%%%%%%%%%%%%%%%%%%%%%%%%%%%%%%%%%%%%%%%%%%%%%%%%%%%%%%%%%%%%%%%%%%%%%%%%
%%%%%%%%%%%%%%%%%%%%%%%%%%%%%%%%%%%%%%%%%%%%%%%%%%%%%%%%%%%%%%%%%%%%%%%%%%%%%%%%%%%%%%%%%%%%%%%%%%%%%%%%%%%%%%%%%%%%%%%%

Let $T$ be a tree.
For an edge $x_{1}x_{2}$ of $T$, let $T_{x_{1}x_{2}}^{x_{1}}$ be the component of $T-x_{1}x_{2}$ containing $x_{1}$.
An edge $x_{1}x_{2}$ of $T$ is {\it good} if $|V(T_{x_{1}x_{2}}^{x_{i}})|\geq 3$ for each $i\in \{1,2\}$.

\begin{figure}
\begin{center}
%WinTpicVersion4.32a
{\unitlength 0.1in%
\begin{picture}(19.0000,9.4500)(3.5000,-12.7500)%
% CIRCLE 2 0 0 0 Black Black  
% 4 1500 1195 1500 1245 1500 1245 1500 1245
% 
\special{sh 1.000}%
\special{ia 1500 1195 50 50 0.0000000 6.2831853}%
\special{pn 8}%
\special{ar 1500 1195 50 50 0.0000000 6.2831853}%
% CIRCLE 2 0 0 0 Black Black  
% 4 800 805 800 855 800 855 800 855
% 
\special{sh 1.000}%
\special{ia 800 805 50 50 0.0000000 6.2831853}%
\special{pn 8}%
\special{ar 800 805 50 50 0.0000000 6.2831853}%
% CIRCLE 2 0 0 0 Black Black  
% 4 1200 805 1200 855 1200 855 1200 855
% 
\special{sh 1.000}%
\special{ia 1200 805 50 50 0.0000000 6.2831853}%
\special{pn 8}%
\special{ar 1200 805 50 50 0.0000000 6.2831853}%
% CIRCLE 2 0 0 0 Black Black  
% 4 1800 805 1800 855 1800 855 1800 855
% 
\special{sh 1.000}%
\special{ia 1800 805 50 50 0.0000000 6.2831853}%
\special{pn 8}%
\special{ar 1800 805 50 50 0.0000000 6.2831853}%
% CIRCLE 2 0 0 0 Black Black  
% 4 2200 795 2200 845 2200 845 2200 845
% 
\special{sh 1.000}%
\special{ia 2200 795 50 50 0.0000000 6.2831853}%
\special{pn 8}%
\special{ar 2200 795 50 50 0.0000000 6.2831853}%
% CIRCLE 2 0 0 0 Black Black  
% 4 800 405 800 455 800 455 800 455
% 
\special{sh 1.000}%
\special{ia 800 405 50 50 0.0000000 6.2831853}%
\special{pn 8}%
\special{ar 800 405 50 50 0.0000000 6.2831853}%
% CIRCLE 2 0 0 0 Black Black  
% 4 1200 405 1200 455 1200 455 1200 455
% 
\special{sh 1.000}%
\special{ia 1200 405 50 50 0.0000000 6.2831853}%
\special{pn 8}%
\special{ar 1200 405 50 50 0.0000000 6.2831853}%
% CIRCLE 2 0 0 0 Black Black  
% 4 1200 395 1200 445 1200 445 1200 445
% 
\special{sh 1.000}%
\special{ia 1200 395 50 50 0.0000000 6.2831853}%
\special{pn 8}%
\special{ar 1200 395 50 50 0.0000000 6.2831853}%
% DOT 0 0 3 0 Black Black  
% 4 1000 580 900 580 1100 580 1100 580
% 
\special{pn 4}%
\special{sh 1}%
\special{ar 1000 580 16 16 0 6.2831853}%
\special{sh 1}%
\special{ar 900 580 16 16 0 6.2831853}%
\special{sh 1}%
\special{ar 1100 580 16 16 0 6.2831853}%
\special{sh 1}%
\special{ar 1100 580 16 16 0 6.2831853}%
% DOT 0 0 3 0 Black Black  
% 4 2000 805 1900 805 2100 805 2100 805
% 
\special{pn 4}%
\special{sh 1}%
\special{ar 2000 805 16 16 0 6.2831853}%
\special{sh 1}%
\special{ar 1900 805 16 16 0 6.2831853}%
\special{sh 1}%
\special{ar 2100 805 16 16 0 6.2831853}%
\special{sh 1}%
\special{ar 2100 805 16 16 0 6.2831853}%
% STR 2 0 3 0 Black Black  
% 4 1500 1240 1500 1340 5 0 0 0
% $x$
\put(15.0000,-13.4000){\makebox(0,0){$x$}}%
% STR 2 0 3 0 Black Black  
% 4 600 695 600 795 5 0 0 0
% $y_{1,1}$
\put(6.0000,-7.9500){\makebox(0,0){$y_{1,1}$}}%
% STR 2 0 3 0 Black Black  
% 4 600 295 600 395 5 0 0 0
% $y_{1,2}$
\put(6.0000,-3.9500){\makebox(0,0){$y_{1,2}$}}%
% STR 2 0 3 0 Black Black  
% 4 1400 295 1400 395 5 0 0 0
% $y_{p,2}$
\put(14.0000,-3.9500){\makebox(0,0){$y_{p,2}$}}%
% STR 2 0 3 0 Black Black  
% 4 1400 695 1400 795 5 0 0 0
% $y_{p,1}$
\put(14.0000,-7.9500){\makebox(0,0){$y_{p,1}$}}%
% STR 2 0 3 0 Black Black  
% 4 2350 705 2350 805 5 0 0 0
% $z_{q}$
\put(23.5000,-8.0500){\makebox(0,0){$z_{q}$}}%
% LINE 2 0 3 0 Black Black  
% 8 1500 1195 800 795 800 795 800 395 1200 395 1200 795 1200 795 1500 1205
% 
\special{pn 8}%
\special{pa 1500 1195}%
\special{pa 800 795}%
\special{fp}%
\special{pa 800 795}%
\special{pa 800 395}%
\special{fp}%
\special{pa 1200 395}%
\special{pa 1200 795}%
\special{fp}%
\special{pa 1200 795}%
\special{pa 1500 1205}%
\special{fp}%
% STR 2 0 3 0 Black Black  
% 4 1660 700 1660 800 5 0 0 0
% $z_{1}$
\put(16.6000,-8.0000){\makebox(0,0){$z_{1}$}}%
% LINE 2 0 3 0 Black Black  
% 4 1800 805 1500 1205 1500 1205 2200 805
% 
\special{pn 8}%
\special{pa 1800 805}%
\special{pa 1500 1205}%
\special{fp}%
\special{pa 1500 1205}%
\special{pa 2200 805}%
\special{fp}%
\end{picture}}%

\caption{Tree $T_{p,q}$}
\label{f-nogood}
\end{center}
\end{figure}
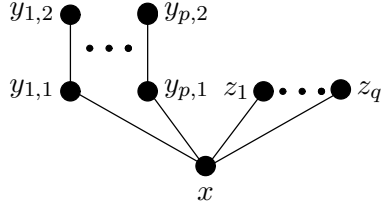

For non-negative integers $p$ and $q$, we let $T_{p,q}$ denote the tree with
$$
V(T_{p,q})=\{x\}\cup \{y_{i,j}:1\leq i\leq p,j\in \{1,2\}\}\cup \{z_{i}:1\leq i\leq q\}
$$
and
$$
E(T_{p,q})=\{xy_{i,1},y_{i,1}y_{i,2}:1\leq i\leq p\}\cup \{xz_{i}:1\leq i\leq q\}
$$
(see Figure~\ref{f-nogood}).

The following lemma might be known.
However, to keep the paper self-contained, we give its proof.

\begin{lem}%%%%%%%%%%%%%%%%%%%%%%%%%%%%%%%%%%%%%%%%%%%%%%%%%%%%%%%%%%%%%%%%%%%%%%%%%%%%%%%%%%%%%%%%%%%%%%%%%%%%%%%%%%%%%
\label{lem2}
Let $T$ be a tree of order at least $3$.
Then $T$ has no good edge if and only if $T$ is isomorphic to $T_{p,q}$ for some $p\geq 0$ and $q\geq 0$ with $2p+q\geq 2$.
\end{lem}
%%%%%%%%%%%%%%%%%%%%%%%%%%%%%%%%%%%%%%%%%%%%%%%%%%%%%%%%%%%%%%%%%%%%%%%%%%%%%%%%%%%%%%%%%%%%%%%%%%%%%%%%%%%%%%%%%%%%%%%%
\proof
For integers $p\geq 0$ and $q\geq 0$ with $2p+q\geq 2$, it is clear that $T_{p,q}$ has no good edge.
Thus it suffices to show that the ``only if'' part of the lemma.

Let $P=x_{0}x_{1}\cdots x_{d}$ be a longest path of $T$.
Then $d$ is equal to the diameter of $T$.
In particular, $d_{T}(x_{1})\geq 2$ because $d\geq 2$.
By the maximality of $P$, every vertex in $N_{T}(x_{1})-\{x_{0},x_{2}\}$ is a leaf of $T$.

Suppose that $d_{T}(x_{1})\geq 3$.
Since $x_{1}x_{2}$ is not a good edge of $T$, $|V(T_{x_{1}x_{2}}^{x_{2}})|\leq 2$.
In particular, either $d=2$ and $V(T_{x_{1}x_{2}}^{x_{2}})=\{x_{2}\}$ or $d=3$ and $V(T_{x_{1}x_{2}}^{x_{2}})=\{x_{2},x_{3}\}$.
Let $k=d_{T}(x_{1})$.
If $d=2$ and $V(T_{x_{1}x_{2}}^{x_{2}})=\{x_{2}\}$, then $T$ is isomorphic to $T_{0,k}$; if $d=3$ and $V(T_{x_{1}x_{2}}^{x_{2}})=\{x_{2},x_{3}\}$, then $T$ is isomorphic to $T_{1,k-1}$.
In either case, we obtain the desired conclusion.
Thus we may assume that $d_{T}(x_{1})=2$ (i.e., $N_{T}(x_{1})=\{x_{0},x_{2}\}$).

For each vertex $u\in N_{T}(x_{2})-\{x_{1}\}$, since $x_{2}u$ is not a good edge and $|V(T_{x_{2}u}^{x_{2}})|\geq 3$, $|V(T_{x_{2}u}^{u})|\leq 2$.
Let $p=|\{u\in N_{T}(x_{2})-\{x_{1}\}:|V(T_{x_{2}u}^{u})|=2\}|$ and $q=|\{u\in N_{T}(x_{2})-\{x_{1}\}:|V(T_{x_{2}u}^{u})|=1\}|$.
Then $T$ is isomorphic to $T_{p+1,q}$, as desired.
\qed

%%%%%%%%%%%%%%%%%%%%%%%%%%%%%%%%%%%%%%%%%%%%%%%%%%%%%%%%%%%%%%%%%%%%%%%%%%%%%%%%%%%%%%%%%%%%%%%%%%%%%%%%%%%%%%%%%%%%%%%%
%%%%%%%%%%%%%%%%%%%%%%%%%%%%%%%%%%%%%%%%%%%%%%%%%%%%%%%%%%%%%%%%%%%%%%%%%%%%%%%%%%%%%%%%%%%%%%%%%%%%%%%%%%%%%%%%%%%%%%%%
%%%%%%%%%%%%%%%%%%%%%%%%%%%%%%%%%%%%%%%%%%%%%%%%%%%%%%%%%%%%%%%%%%%%%%%%%%%%%%%%%%%%%%%%%%%%%%%%%%%%%%%%%%%%%%%%%%%%%%%%
\section{Upper bound on $\gamma ^{c}$ for $\frac{1}{2}\leq c<1$}\label{sec3}
%%%%%%%%%%%%%%%%%%%%%%%%%%%%%%%%%%%%%%%%%%%%%%%%%%%%%%%%%%%%%%%%%%%%%%%%%%%%%%%%%%%%%%%%%%%%%%%%%%%%%%%%%%%%%%%%%%%%%%%%
%%%%%%%%%%%%%%%%%%%%%%%%%%%%%%%%%%%%%%%%%%%%%%%%%%%%%%%%%%%%%%%%%%%%%%%%%%%%%%%%%%%%%%%%%%%%%%%%%%%%%%%%%%%%%%%%%%%%%%%%
%%%%%%%%%%%%%%%%%%%%%%%%%%%%%%%%%%%%%%%%%%%%%%%%%%%%%%%%%%%%%%%%%%%%%%%%%%%%%%%%%%%%%%%%%%%%%%%%%%%%%%%%%%%%%%%%%%%%%%%%

In this section, we prove the following theorem.

\begin{thm}%%%%%%%%%%%%%%%%%%%%%%%%%%%%%%%%%%%%%%%%%%%%%%%%%%%%%%%%%%%%%%%%%%%%%%%%%%%%%%%%%%%%%%%%%%%%%%%%%%%%%%%%%%%%%
\label{thm3.1}
Let $c\in \mathbb{R}^{+}$ be a number with $\frac{1}{2}\leq c<1$, and let $m\geq 1$ be the integer such that $\frac{m}{m+1}\leq c<\frac{m+1}{m+2}$.
Let $T$ be a tree of order $n\geq 3$.
Then the conclusion of Theorem~\ref{mainthm} holds for $G=T$.
%Let $T$ be a tree of order $n\geq 3$.
%Then the following hold:
%\begin{enumerate}
%\item[{\upshape(i)}]
%If $\frac{m}{m+1}\leq c<\frac{2m+1}{2m+3}$, then $\gamma ^{c}(T)\leq \frac{m+1}{2m+3}n$; and
%$\item[{\upshape(ii)}]
%if $\frac{2m+1}{2m+3}\leq c<\frac{m+1}{m+2}$, then $\gamma ^{c}(T)\leq \frac{cm+2c+1}{2m+5}n$.
%\end{enumerate}
\end{thm}
%%%%%%%%%%%%%%%%%%%%%%%%%%%%%%%%%%%%%%%%%%%%%%%%%%%%%%%%%%%%%%%%%%%%%%%%%%%%%%%%%%%%%%%%%%%%%%%%%%%%%%%%%%%%%%%%%%%%%%%%

\begin{remark}%%%%%%%%%%%%%%%%%%%%%%%%%%%%%%%%%%%%%%%%%%%%%%%%%%%%%%%%%%%%%%%%%%%%%%%%%%%%%%%%%%%%%%%%%%%%%%%%%%%%%%%%%%
Since the deletion of edges cannot decrease the $c$-self-domination number of a graph, we obtain Theorem~\ref{mainthm} for the case where $\frac{1}{2}\leq c<1$ as a corollary of Theorem~\ref{thm3.1}.
\end{remark}
%%%%%%%%%%%%%%%%%%%%%%%%%%%%%%%%%%%%%%%%%%%%%%%%%%%%%%%%%%%%%%%%%%%%%%%%%%%%%%%%%%%%%%%%%%%%%%%%%%%%%%%%%%%%%%%%%%%%%%%%

We first show the following lemma.

\begin{lem}%%%%%%%%%%%%%%%%%%%%%%%%%%%%%%%%%%%%%%%%%%%%%%%%%%%%%%%%%%%%%%%%%%%%%%%%%%%%%%%%%%%%%%%%%%%%%%%%%%%%%%%%%%%%%
\label{lem3.1}
Let $c$ and $m$ be as in Theorem~\ref{thm3.1}.
If $T=T_{p,q}$ for integers $p\geq 0$ and $q\geq 0$ with $2p+q\geq 2$, then the conclusion of Theorem~\ref{thm3.1} holds.
\end{lem}
%%%%%%%%%%%%%%%%%%%%%%%%%%%%%%%%%%%%%%%%%%%%%%%%%%%%%%%%%%%%%%%%%%%%%%%%%%%%%%%%%%%%%%%%%%%%%%%%%%%%%%%%%%%%%%%%%%%%%%%%
\proof
The function $f_{1}:V(T_{p,q})\rightarrow \mathbb{R}^{+}$ with
$$
f_{1}(u)=
\begin{cases}
1, & u=x;\\
c, & u\in \{y_{i,2}:1\leq i\leq p\};\\
0, & \mbox{otherwise}
\end{cases}
$$
is a $c$-SDF of $T_{p,q}$ with $w(f_{1})=cp+1$.
Hence
\begin{align}
\gamma ^{c}(T_{p,q})\leq w(f_{1})=cp+1=\frac{cp+1}{2p+q+1}|V(T_{p,q})|.\label{eq3.1-1}
\end{align}
If $p\geq 1$, then the function $f_{2}:V(T_{p,q})\rightarrow \mathbb{R}^{+}$ with
$$
f_{2}(u)=
\begin{cases}
1, & u\in \{y_{i,1}:1\leq i\leq p\};\\
c, & u\in \{z_{i}:1\leq i\leq q\};\\
0, & \mbox{otherwise}
\end{cases}
$$
is a $c$-SDF of $T_{p,q}$ with $w(f_{2})=p+cq$.
Hence
\begin{align}
\gamma ^{c}(T_{p,q})\leq w(f_{2})=p+cq=\frac{p+cq}{2p+q+1}|V(T_{p,q})|\mbox{~~~~if }p\geq 1.\label{eq3.1-2}
\end{align}
Furthermore, we have
\begin{align}
\frac{m+1}{2m+3}=\frac{\frac{(2m+1)(m+2)}{2m+3}+1}{2m+5}\leq \frac{c(m+2)+1}{2m+5}\mbox{~~~~if }c\geq \frac{2m+1}{2m+3}.\label{eq3.1-3}
\end{align}

We divide the proof into four cases.

\medskip
\noindent
\textbf{Case 1:} Either $p=0$ or $(p,q)=(1,0)$.

Note that if $p=0$, then $q=2p+q\geq 2$.
Hence it follows from (\ref{eq3.1-1}) and (\ref{eq3.1-2}) that
\begin{align}
\gamma ^{c}(T_{p,q}) \leq \frac{1}{3}|V(T_{p,q})|<\frac{m+1}{2m+3}|V(T_{p,q})|.\label{eq3.1-4}
\end{align}
In particular, we obtain the desired conclusion for the case where $\frac{m}{m+1}\leq c<\frac{2m+1}{2m+3}$.
If $\frac{2m+1}{2m+3}\leq c<\frac{m+1}{m+2}$, then it follows from (\ref{eq3.1-3}) and (\ref{eq3.1-4}) that
$$
\gamma ^{c}(T_{p,q})<\frac{m+1}{2m+3}|V(T_{p,q})|\leq \frac{c(m+2)+1}{2m+5}|V(T_{p,q})|,
$$
which follows the desired conclusion for the remaining case.

\medskip
\noindent
\textbf{Case 2:} $p=1$ and $q\geq 1$.

If $\frac{m}{m+1}\leq c<\frac{2m+1}{2m+3}$, then by (\ref{eq3.1-1}) and the assumption that $q\geq 1$,
$$
\gamma ^{c}(T_{p,q})\leq \frac{cp+1}{2p+q+1}|V(T_{p,q})|<\frac{\frac{2m+1}{2m+3}+1}{4}|V(T_{p,q})|=\frac{m+1}{2m+3}|V(T_{p,q})|,
$$
as desired.

Next we suppose that $\frac{2m+1}{2m+3}\leq c<\frac{m+1}{m+2}$.
Since
\begin{align*}
4(cm+2c+1)-(c+1)(2m+5) &= c(2m+3)-(2m+1)\\
&\geq \frac{(2m+1)(2m+3)}{2m+3}-(2m+1)\\
&=0,
\end{align*}
we have
$$
\frac{cp+1}{2p+q+1}\leq \frac{c+1}{4}\leq \frac{cm+2c+1}{2m+5}.
$$
This together with (\ref{eq3.1-1}) leads to the desired conclusion for the case where $\frac{2m+1}{2m+3}\leq c<\frac{m+1}{m+2}$.

\medskip
\noindent
\textbf{Case 3:} $p\geq 2$ and $q\geq 1$.

Since $p\geq 2$,
$$
2(m+1)(m+2)(p+1)-((m+1)p+m+2)(2m+3)=m(p-1)+p-2>0,
$$
and hence
$$
\frac{(m+1)p+m+2}{2(m+2)(p+1)}<\frac{m+1}{2m+3}.
$$
This together with (\ref{eq3.1-1}) implies
\begin{align}
\gamma ^{c}(T_{p,q}) &\leq \frac{cp+1}{2p+q+1}|V(T_{p,q})|\nonumber \\
&< \frac{\frac{(m+1)p}{m+2}+1}{2p+2}|V(T_{p,q})|\nonumber \\
&= \frac{(m+1)p+m+2}{2(m+2)(p+1)}|V(T_{p,q})|\nonumber \\
&< \frac{m+1}{2m+3}|V(T_{p,q})|.\label{eq3.1-5}
\end{align}
In particular, we obtain the desired conclusion for the case where $\frac{m}{m+1}\leq c<\frac{2m+1}{2m+3}$.
If $\frac{2m+1}{2m+3}\leq c<\frac{m+1}{m+2}$, then it follows from (\ref{eq3.1-3}) and (\ref{eq3.1-5}) that
$$
\gamma ^{c}(T_{p,q})\leq \frac{m+1}{2m+3}|V(T_{p,q})|\leq \frac{c(m+2)+1}{2m+5}|V(T_{p,q})|,
$$
which follows the desired conclusion for the remaining case.

\medskip
\noindent
\textbf{Case 4:} $p\geq 2$ and $q=0$.

Suppose that $p\leq m+1$.
Then it follows from (\ref{eq3.1-2}) that
\begin{align}
\gamma ^{c}(T_{p,q})\leq \frac{p}{2p+1}|V(T_{p,q})|\leq \frac{m+1}{2m+3}|V(T_{p,q})|.\label{eq3.1-6}
\end{align}
In particular, we obtain the desired conclusion for the case where $\frac{m}{m+1}\leq c<\frac{2m+1}{2m+3}$.
If $\frac{2m+1}{2m+3}\leq c<\frac{m+1}{m+2}$, then it follows from (\ref{eq3.1-3}) and (\ref{eq3.1-6}) that
$$
\gamma ^{c}(T_{p,q})\leq \frac{m+1}{2m+3}|V(T_{p,q})|\leq \frac{c(m+2)+1}{2m+5}|V(T_{p,q})|,
$$
as desired.
Thus we may assume that $p\geq m+2$.

We first suppose that $\frac{m}{m+1}\leq c<\frac{2m+1}{2m+3}$.
Then
$$
(m+1)(2p+1)-((2m+1)p+2m+3)=p-m-2\geq 0,
$$
and hence
$$
\frac{(2m+1)p+2m+3}{2p+1}\leq m+1.
$$
This together with (\ref{eq3.1-1}) leads to
\begin{align*}
\gamma ^{c}(T_{p,q}) &\leq \frac{cp+1}{2p+1}|V(T_{p,q})|\\
&< \frac{\frac{(2m+1)p}{2m+3}+1}{2p+1}|V(T_{p,q})|\\
&= \frac{(2m+1)p+2m+3}{(2p+1)(2m+3)}|V(T_{p,q})|\\
&\leq \frac{m+1}{2m+3}|V(T_{p,q})|,
\end{align*}
which leads to the desired conclusion for the case where $\frac{m}{m+1}\leq c<\frac{2m+1}{2m+3}$.

Next we suppose that $\frac{2m+1}{2m+3}\leq c<\frac{m+1}{m+2}$.
Then
$$
(c(m+2)+1)(2p+1)-(cp+1)(2m+5)=(p-m-2)(2-c)\geq 0,
$$
and hence
$$
\frac{cp+1}{2p+1}\leq \frac{c(m+2)+1}{2m+5}.
$$
This together with (\ref{eq3.1-1}) leads to
$$
\gamma ^{c}(T_{p,q})\leq \frac{cp+1}{2p+1}|V(T_{p,q})|\leq \frac{c(m+2)+1}{2m+5}|V(T_{p,q})|,
$$
which leads to the desired conclusion for the case where $\frac{2m+1}{2m+3}\leq c<\frac{m+1}{m+2}$.

This completes the proof of Lemma~\ref{lem3.1}.
\qed

\medbreak\noindent\textit{Proof of Theorem~\ref{thm3.1}.}\quad
We proceed by induction on $n$.
If $T$ has no good edge, then by Lemma~\ref{lem2} and Lemma~\ref{lem3.1}, the desired conclusion holds.
Thus we may assume that $T$ has a good edge $x_{1}x_{2}$.
Then by the induction hypothesis, for each $i\in \{1,2\}$,
$$
\gamma ^{c}(T_{x_{1}x_{2}}^{x_{i}})\leq
\begin{cases}
\frac{m+1}{2m+3}|V(T_{x_{1}x_{2}}^{x_{i}})|, & \frac{m}{m+1}\leq c<\frac{2m+1}{2m+3};\\
\frac{cm+2c+1}{2m+5}|V(T_{x_{1}x_{2}}^{x_{i}})|, & \frac{2m+1}{2m+3}\leq c<\frac{m+1}{m+2}.
\end{cases}
$$
Since $\gamma ^{c}(T)\leq \gamma ^{c}(T_{x_{1}x_{2}}^{x_{1}})+\gamma ^{c}(T_{x_{1}x_{2}}^{x_{2}})$ and $|V(T_{x_{1}x_{2}}^{x_{1}})|+|V(T_{x_{1}x_{2}}^{x_{2}})|=n$, this leads to the desired conclusion.
\qed

%%%%%%%%%%%%%%%%%%%%%%%%%%%%%%%%%%%%%%%%%%%%%%%%%%%%%%%%%%%%%%%%%%%%%%%%%%%%%%%%%%%%%%%%%%%%%%%%%%%%%%%%%%%%%%%%%%%%%%%%
%%%%%%%%%%%%%%%%%%%%%%%%%%%%%%%%%%%%%%%%%%%%%%%%%%%%%%%%%%%%%%%%%%%%%%%%%%%%%%%%%%%%%%%%%%%%%%%%%%%%%%%%%%%%%%%%%%%%%%%%
%%%%%%%%%%%%%%%%%%%%%%%%%%%%%%%%%%%%%%%%%%%%%%%%%%%%%%%%%%%%%%%%%%%%%%%%%%%%%%%%%%%%%%%%%%%%%%%%%%%%%%%%%%%%%%%%%%%%%%%%
\section{Upper bound on $\gamma ^{c}$ for $1<c\leq 2$}\label{sec4}
%%%%%%%%%%%%%%%%%%%%%%%%%%%%%%%%%%%%%%%%%%%%%%%%%%%%%%%%%%%%%%%%%%%%%%%%%%%%%%%%%%%%%%%%%%%%%%%%%%%%%%%%%%%%%%%%%%%%%%%%
%%%%%%%%%%%%%%%%%%%%%%%%%%%%%%%%%%%%%%%%%%%%%%%%%%%%%%%%%%%%%%%%%%%%%%%%%%%%%%%%%%%%%%%%%%%%%%%%%%%%%%%%%%%%%%%%%%%%%%%%
%%%%%%%%%%%%%%%%%%%%%%%%%%%%%%%%%%%%%%%%%%%%%%%%%%%%%%%%%%%%%%%%%%%%%%%%%%%%%%%%%%%%%%%%%%%%%%%%%%%%%%%%%%%%%%%%%%%%%%%%

In this section, we prove the following theorem.

\begin{thm}%%%%%%%%%%%%%%%%%%%%%%%%%%%%%%%%%%%%%%%%%%%%%%%%%%%%%%%%%%%%%%%%%%%%%%%%%%%%%%%%%%%%%%%%%%%%%%%%%%%%%%%%%%%%%
\label{thm4.1}
Let $c\in \mathbb{R}^{+}$ be a number with $1<c\leq 2$, and let $m\geq 1$ be the integer such that $\frac{m+2}{m+1}<c\leq \frac{m+1}{m}$.
Let $T$ be a tree of order $n\geq 3$.
Then the conclusion of Theorem~\ref{mainthm} holds for $G=T$.
%Then the following hold:
%\begin{enumerate}[{\upshape(i)}]
%\item
%If $\frac{m+2}{m+1}<c\leq \frac{(2m+1)(m+2)}{(2m+3)m}$, then $\gamma ^{c}(T)\leq \frac{m+2}{2m+3}n$; and
%\item
%if $\frac{(2m+1)(m+2)}{(2m+3)m}<c\leq \frac{m+1}{m}$, then $\gamma ^{c}(T)\leq \frac{cm}{2m+1}n$.
%\end{enumerate}
\end{thm}
%%%%%%%%%%%%%%%%%%%%%%%%%%%%%%%%%%%%%%%%%%%%%%%%%%%%%%%%%%%%%%%%%%%%%%%%%%%%%%%%%%%%%%%%%%%%%%%%%%%%%%%%%%%%%%%%%%%%%%%%

\begin{remark}%%%%%%%%%%%%%%%%%%%%%%%%%%%%%%%%%%%%%%%%%%%%%%%%%%%%%%%%%%%%%%%%%%%%%%%%%%%%%%%%%%%%%%%%%%%%%%%%%%%%%%%%%%
Since the deletion of edges cannot decrease the $c$-self-domination number of a graph, we obtain Theorem~\ref{mainthm} for the case where $1<c\leq 2$ as a corollary of Theorem~\ref{thm4.1}.
\end{remark}
%%%%%%%%%%%%%%%%%%%%%%%%%%%%%%%%%%%%%%%%%%%%%%%%%%%%%%%%%%%%%%%%%%%%%%%%%%%%%%%%%%%%%%%%%%%%%%%%%%%%%%%%%%%%%%%%%%%%%%%%

We first show the following lemma.

\begin{lem}%%%%%%%%%%%%%%%%%%%%%%%%%%%%%%%%%%%%%%%%%%%%%%%%%%%%%%%%%%%%%%%%%%%%%%%%%%%%%%%%%%%%%%%%%%%%%%%%%%%%%%%%%%%%%
\label{lem4.1}
Let $c$ and $m$ be as in Theorem~\ref{thm4.1}.
If $T=T_{p,q}$ for integers $p\geq 0$ and $q\geq 0$ with $2p+q\geq 2$, then the conclusion of Theorem~\ref{thm4.1} holds.
\end{lem}
%%%%%%%%%%%%%%%%%%%%%%%%%%%%%%%%%%%%%%%%%%%%%%%%%%%%%%%%%%%%%%%%%%%%%%%%%%%%%%%%%%%%%%%%%%%%%%%%%%%%%%%%%%%%%%%%%%%%%%%%
\proof
If $p=0$, then the function $f_{1}:$
The function $f_{1}:V(T_{p,q})\rightarrow \mathbb{R}^{+}$ with
$$
f_{1}(u)=
\begin{cases}
c, & u=x;\\
0, & \mbox{otherwise}
\end{cases}
$$
is a $c$-SDF of $T_{p,q}$ with $w(f_{1})=c$.
Hence
\begin{align}
\gamma ^{c}(T_{p,q})\leq w(f_{1})=c=\frac{c}{q+1}|V(T_{p,q})|\mbox{~~~~if }p=0.\label{eq4.1-1}
\end{align}
If $p\geq 1$, then the function $f_{2}:V(T_{p,q})\rightarrow \mathbb{R}^{+}$ with
$$
f_{2}(u)=
\begin{cases}
1, & u\in \{x,y_{i,1}:1\leq i\leq p\};\\
0, & \mbox{otherwise}
\end{cases}
$$
is a $c$-SDF of $T_{p,q}$ with $w(f_{2})=p+1$.
Hence
\begin{align}
\gamma ^{c}(T_{p,q})\leq w(f_{2})=p+1=\frac{p+1}{2p+q+1}|V(T_{p,q})|\mbox{~~~~if }p\geq 1.\label{eq4.1-2}
\end{align}
If $q=0$ (i.e., $p\geq 1$), then the function $f_{3}:V(T_{p,q})\rightarrow \mathbb{R}^{+}$ with
$$
f_{3}(u)=
\begin{cases}
c, & u\in \{y_{i,1}:1\leq i\leq p\};\\
0, & \mbox{otherwise}
\end{cases}
$$
is a $c$-SDF of $T_{p,q}$ with $w(f_{3})=cp$.
Hence
\begin{align}
\gamma ^{c}(T_{p,q})\leq w(f_{3})=cp=\frac{cp}{2p+1}|V(T_{p,q})|\mbox{~~~~if }q=0.\label{eq4.1-3}
\end{align}
Furthermore, we have
\begin{align}
\frac{cm}{2m+1}\leq \frac{\frac{(2m+1)(m+2)m}{(2m+3)m}}{2m+1}=\frac{m+2}{2m+3}\mbox{~~~~if }c\leq \frac{(2m+1)(m+2)}{(2m+3)m}.\label{eq4.1-4}
\end{align}

We divide the proof into three cases.

\medskip
\noindent
\textbf{Case 1:} Either $p=0$ or $(p,q)=(1,0)$.

Note that if $p=0$, then $q=2p+q\geq 2$.
Hence it follows from (\ref{eq4.1-1}) and (\ref{eq4.1-3}) that
\begin{align}
\gamma ^{c}(T_{p,q}) \leq \frac{c}{3}|V(T_{p,q})|\leq \frac{cm}{2m+1}|V(T_{p,q})|.\label{eq4.1-5}
\end{align}
In particular, we obtain the desired conclusion for the case where $\frac{(2m+1)(m+2)}{(2m+3)m}<c\leq \frac{m+1}{m}$.
If $\frac{m+2}{m+1}<c\leq \frac{(2m+1)(m+2)}{(2m+3)m}$, then it follows from (\ref{eq4.1-4}) and (\ref{eq4.1-5}) that
$$
\gamma ^{c}(T_{p,q})\leq \frac{cm}{2m+1}|V(T_{p,q})|\leq \frac{m+2}{2m+3}|V(T_{p,q})|,
$$
which follows the desired conclusion for the remaining case.

\medskip
\noindent
\textbf{Case 2:} $p\geq 1$ and $q\geq 1$.

By (\ref{eq4.1-2}),
\begin{align}
\gamma ^{c}(T_{p,q}) &\leq \frac{p+1}{2p+q+1}|V(T_{p,q})|\nonumber \\
&\leq \frac{p+1}{2p+2}|V(T_{p,q})|\nonumber \\
&= \frac{1}{2}|V(T_{p,q})|\nonumber \\
&\leq \frac{(m+2)m}{(m+1)(2m+1)}|V(T_{p,q})|\nonumber \\
&< \frac{cm}{2m+1}|V(T_{p,q})|.\label{eq4.1-6}
\end{align}
In particular, we obtain the desired conclusion for the case where $\frac{(2m+1)(m+2)}{(2m+3)m}<c\leq \frac{m+1}{m}$.
If $\frac{m+2}{m+1}<c\leq \frac{(2m+1)(m+2)}{(2m+3)m}$, then it follows from (\ref{eq4.1-4}) and (\ref{eq4.1-6}) that
$$
\gamma ^{c}(T_{p,q})<\frac{cm}{2m+1}|V(T_{p,q})|\leq \frac{m+2}{2m+3}|V(T_{p,q})|,
$$
which follows the desired conclusion for the remaining case.

\medskip
\noindent
\textbf{Case 3:} $p\geq 2$ and $q=0$.

Suppose that $p\leq m$.
Then it follows from (\ref{eq4.1-3}) that
\begin{align}
\gamma ^{c}(T_{p,q}) \leq \frac{cp}{2p+1}|V(T_{p,q})|\leq \frac{cm}{2m+1}|V(T_{p,q})|.\label{eq4.1-7}
\end{align}
In particular, we obtain the desired conclusion for the case where $\frac{(2m+1)(m+2)}{(2m+3)m}<c\leq \frac{m+1}{m}$.
If $\frac{m+2}{m+1}<c\leq \frac{(2m+1)(m+2)}{(2m+3)m}$, then it follows from (\ref{eq4.1-4}) and (\ref{eq4.1-7}) that
$$
\gamma ^{c}(T_{p,q})\leq \frac{cm}{2m+1}|V(T_{p,q})|\leq \frac{m+2}{2m+3}|V(T_{p,q})|,
$$
which follows the desired conclusion for the case where $\frac{m+2}{m+1}<c\leq \frac{(2m+1)(m+2)}{(2m+3)m}$.
Thus we may assume that $p\geq m+1$.

We first suppose that $\frac{m+2}{m+1}<c\leq \frac{(2m+1)(m+2)}{(2m+3)m}$.
Then $\frac{p+1}{2p+1}\leq \frac{m+2}{2m+3}$.
This together with (\ref{eq4.1-2}) leads to
$$
\gamma ^{c}(T_{p,q})\leq \frac{p+1}{2p+1}|V(T_{p,q})|\leq \frac{m+2}{2m+3}|V(T_{p,q})|,
$$
as desired.

Next we suppose that $\frac{(2m+1)(m+2)}{(2m+3)m}<c\leq \frac{m+1}{m}$.
Then
$$
\frac{p+1}{2p+1}\leq \frac{m+2}{2m+3}=\frac{\frac{(2m+1)(m+2)m}{(2m+3)m}}{2m+1}<\frac{cm}{2m+1}.
$$
This together with (\ref{eq4.1-2}) leads to
$$
\gamma ^{c}(T_{p,q}) \leq \frac{p+1}{2p+1}|V(T_{p,q})|<\frac{cm}{2m+1}|V(T_{p,q})|,
$$
as desired.

This completes the proof of Lemma~\ref{lem4.1}.
\qed

\medbreak\noindent\textit{Proof of Theorem~\ref{thm4.1}.}\quad
We proceed by induction on $n$.
If $T$ has no good edge, then by Lemma~\ref{lem2} and Lemma~\ref{lem4.1}, the desired conclusion holds.
Thus we may assume that $T$ has a good edge $x_{1}x_{2}$.
Then by the induction hypothesis, for each $i\in \{1,2\}$,
$$
\gamma ^{c}(T_{x_{1}x_{2}}^{x_{i}})\leq
\begin{cases}
\frac{m+2}{2m+3}|V(T_{x_{1}x_{2}}^{x_{i}})|, & \frac{m+2}{m+1}<c\leq \frac{(2m+1)(m+2)}{(2m+3)m};\\
\frac{cm}{2m+1}|V(T_{x_{1}x_{2}}^{x_{i}})|, & \frac{(2m+1)(m+2)}{(2m+3)m}<c\leq \frac{m+1}{m}.
\end{cases}
$$
Since $\gamma ^{c}(T)\leq \gamma ^{c}(T_{x_{1}x_{2}}^{x_{1}})+\gamma ^{c}(T_{x_{1}x_{2}}^{x_{2}})$ and $|V(T_{x_{1}x_{2}}^{x_{1}})|+|V(T_{x_{1}x_{2}}^{x_{2}})|=n$, this leads to the desired conclusion.
\qed

%%%%%%%%%%%%%%%%%%%%%%%%%%%%%%%%%%%%%%%%%%%%%%%%%%%%%%%%%%%%%%%%%%%%%%%%%%%%%%%%%%%%%%%%%%%%%%%%%%%%%%%%%%%%%%%%%%%%%%%%
%%%%%%%%%%%%%%%%%%%%%%%%%%%%%%%%%%%%%%%%%%%%%%%%%%%%%%%%%%%%%%%%%%%%%%%%%%%%%%%%%%%%%%%%%%%%%%%%%%%%%%%%%%%%%%%%%%%%%%%%
%%%%%%%%%%%%%%%%%%%%%%%%%%%%%%%%%%%%%%%%%%%%%%%%%%%%%%%%%%%%%%%%%%%%%%%%%%%%%%%%%%%%%%%%%%%%%%%%%%%%%%%%%%%%%%%%%%%%%%%%
\section{Upper bound on $\gamma ^{c}$ for $c>2$}\label{sec5}
%%%%%%%%%%%%%%%%%%%%%%%%%%%%%%%%%%%%%%%%%%%%%%%%%%%%%%%%%%%%%%%%%%%%%%%%%%%%%%%%%%%%%%%%%%%%%%%%%%%%%%%%%%%%%%%%%%%%%%%%
%%%%%%%%%%%%%%%%%%%%%%%%%%%%%%%%%%%%%%%%%%%%%%%%%%%%%%%%%%%%%%%%%%%%%%%%%%%%%%%%%%%%%%%%%%%%%%%%%%%%%%%%%%%%%%%%%%%%%%%%
%%%%%%%%%%%%%%%%%%%%%%%%%%%%%%%%%%%%%%%%%%%%%%%%%%%%%%%%%%%%%%%%%%%%%%%%%%%%%%%%%%%%%%%%%%%%%%%%%%%%%%%%%%%%%%%%%%%%%%%%

\begin{lem}%%%%%%%%%%%%%%%%%%%%%%%%%%%%%%%%%%%%%%%%%%%%%%%%%%%%%%%%%%%%%%%%%%%%%%%%%%%%%%%%%%%%%%%%%%%%%%%%%%%%%%%%%%%%%
\label{lem5.1}
Let $c\in \mathbb{R}^{+}$ be a number with $c\geq 2$.
Let $G$ be a connected graph of order at least $2$.
Then $\gamma ^{c}(G)=\gamma _{t}(G)$.
\end{lem}
%%%%%%%%%%%%%%%%%%%%%%%%%%%%%%%%%%%%%%%%%%%%%%%%%%%%%%%%%%%%%%%%%%%%%%%%%%%%%%%%%%%%%%%%%%%%%%%%%%%%%%%%%%%%%%%%%%%%%%%%
\proof
Let $f$ be a $\gamma ^{c}$-function of $G$.
By (\ref{01c}), we have $\{f(u):u\in V(G)\}\subseteq \{0,1,c\}$.
Choose $f$ so that $|\{u\in V(G):f(u)=1\}|$ is as large as possible.

Suppose that $f(x)=c$ for some $x\in V(G)$, and let $y\in N_{G}(x)$.
Then the function $f':V(G)\rightarrow \mathbb{R}^{+}$ with
$$
f'(u)=
\begin{cases}
1, & u\in \{x,y\};\\
f(u), & \mbox{otherwise}
\end{cases}
$$
is a $c$-SDF of $G$ with $w(f')=w(f)-c-f(y)+2\leq w(f)$ and $|\{u\in V(G):f'(u)=1\}|>|\{u\in V(G):f(u)=1\}|$, which contradicts the choice of $f$.
Thus $\{f(u):u\in V(G)\}\subseteq \{0,1\}$.

Since the set $S:=\{u\in V(G):f(u)=1\}$ is a total dominating set of $G$ with $|S|=w(f)$, $\gamma _{t}(G)\leq |S|=w(f)=\gamma ^{c}(G)$.
Since every $\infty $-SDF of $G$ is also a $c$-SDF of $G$, we have $\gamma ^{c}(G)\leq \gamma ^{\infty }(G)$.
This together with Proposition~\ref{prop1.1} leads to $\gamma ^{c}(G)\leq \gamma ^{\infty }(G)=\gamma _{t}(G)$, we obtain the desired conclusion.
\qed

By Lemma~\ref{lem5.1} and Theorem~\ref{thm4.1} for the case where $c\geq 2$, we obtain the following proposition.

\begin{prop}%%%%%%%%%%%%%%%%%%%%%%%%%%%%%%%%%%%%%%%%%%%%%%%%%%%%%%%%%%%%%%%%%%%%%%%%%%%%%%%%%%%%%%%%%%%%%%%%%%%%%%%%%%%%
\label{prop5.1}
Let $c\in \mathbb{R}^{+}$ be a number with $c\geq 2$.
Let $G$ be a connected graph of order $n\geq 3$.
Then $\gamma ^{c}(G)\leq \frac{2}{3}n$.
\end{prop}
%%%%%%%%%%%%%%%%%%%%%%%%%%%%%%%%%%%%%%%%%%%%%%%%%%%%%%%%%%%%%%%%%%%%%%%%%%%%%%%%%%%%%%%%%%%%%%%%%%%%%%%%%%%%%%%%%%%%%%%%

\begin{remark}%%%%%%%%%%%%%%%%%%%%%%%%%%%%%%%%%%%%%%%%%%%%%%%%%%%%%%%%%%%%%%%%%%%%%%%%%%%%%%%%%%%%%%%%%%%%%%%%%%%%%%%%%%
We obtain Theorem~\ref{mainthm} for the case where $c>2$ as a corollary of Proposition~\ref{prop5.1}.
\end{remark}
%%%%%%%%%%%%%%%%%%%%%%%%%%%%%%%%%%%%%%%%%%%%%%%%%%%%%%%%%%%%%%%%%%%%%%%%%%%%%%%%%%%%%%%%%%%%%%%%%%%%%%%%%%%%%%%%%%%%%%%%

%%%%%%%%%%%%%%%%%%%%%%%%%%%%%%%%%%%%%%%%%%%%%%%%%%%%%%%%%%%%%%%%%%%%%%%%%%%%%%%%%%%%%%%%%%%%%%%%%%%%%%%%%%%%%%%%%%%%%%%%
%%%%%%%%%%%%%%%%%%%%%%%%%%%%%%%%%%%%%%%%%%%%%%%%%%%%%%%%%%%%%%%%%%%%%%%%%%%%%%%%%%%%%%%%%%%%%%%%%%%%%%%%%%%%%%%%%%%%%%%%
%%%%%%%%%%%%%%%%%%%%%%%%%%%%%%%%%%%%%%%%%%%%%%%%%%%%%%%%%%%%%%%%%%%%%%%%%%%%%%%%%%%%%%%%%%%%%%%%%%%%%%%%%%%%%%%%%%%%%%%%
\section{Examples}\label{sec6}
%%%%%%%%%%%%%%%%%%%%%%%%%%%%%%%%%%%%%%%%%%%%%%%%%%%%%%%%%%%%%%%%%%%%%%%%%%%%%%%%%%%%%%%%%%%%%%%%%%%%%%%%%%%%%%%%%%%%%%%%
%%%%%%%%%%%%%%%%%%%%%%%%%%%%%%%%%%%%%%%%%%%%%%%%%%%%%%%%%%%%%%%%%%%%%%%%%%%%%%%%%%%%%%%%%%%%%%%%%%%%%%%%%%%%%%%%%%%%%%%%
%%%%%%%%%%%%%%%%%%%%%%%%%%%%%%%%%%%%%%%%%%%%%%%%%%%%%%%%%%%%%%%%%%%%%%%%%%%%%%%%%%%%%%%%%%%%%%%%%%%%%%%%%%%%%%%%%%%%%%%%

In this section, we show that Theorem~\ref{mainthm} is best possible.

Let $p\geq 1$ and $s\geq 1$ be integers, and let $T_{p,0}$ be the tree defined in Section~\ref{sec2}.
Let $L_{p}^{1},\ldots ,L_{p}^{s}$ be vertex-disjoint copies of $T_{p,0}$.
For $l~(1\leq l\leq s)$, let $x^{l}$ and $y_{i,j}^{l}~(1\leq i\leq p,j\in \{1,2\})$ be the vertices of $L_{p}^{l}$ corresponding to $x$ and $y_{i,j}$, respectively.
Let $T_{p}^{s}$ be the tree obtained from $L_{p}^{1},\ldots ,L_{p}^{s}$ by adding edges $x^{l}x^{l+1}~(1\leq l\leq s-1)$.

%%%%%%%%%%%%%%%%%%%%%%%%%%%%%%%%%%%%%%%%%%%%%%%%%%%%%%%%%%%%%%%%%%%%%%%%%%%%%%%%%%%%%%%%%%%%%%%%%%%%%%%%%%%%%%%%%%%%%%%%
%%%%%%%%%%%%%%%%%%%%%%%%%%%%%%%%%%%%%%%%%%%%%%%%%%%%%%%%%%%%%%%%%%%%%%%%%%%%%%%%%%%%%%%%%%%%%%%%%%%%%%%%%%%%%%%%%%%%%%%%
\subsection{The case $\frac{1}{2}\leq c<1$}\label{sec6.1}
%%%%%%%%%%%%%%%%%%%%%%%%%%%%%%%%%%%%%%%%%%%%%%%%%%%%%%%%%%%%%%%%%%%%%%%%%%%%%%%%%%%%%%%%%%%%%%%%%%%%%%%%%%%%%%%%%%%%%%%%
%%%%%%%%%%%%%%%%%%%%%%%%%%%%%%%%%%%%%%%%%%%%%%%%%%%%%%%%%%%%%%%%%%%%%%%%%%%%%%%%%%%%%%%%%%%%%%%%%%%%%%%%%%%%%%%%%%%%%%%%

Throughout this subsection, fix a number $c$ with $\frac{1}{2}\leq c<1$, and let $m\geq 1$ be the integer such that $\frac{m}{m+1}\leq c<\frac{m+1}{m+2}$.

\begin{lem}%%%%%%%%%%%%%%%%%%%%%%%%%%%%%%%%%%%%%%%%%%%%%%%%%%%%%%%%%%%%%%%%%%%%%%%%%%%%%%%%%%%%%%%%%%%%%%%%%%%%%%%%%%%%%
\label{lem6.1-1}
Let $p\geq 1$ and $s\geq 1$ be integers.
Let $f:V(T_{p}^{s})\rightarrow \mathbb{R}^{+}$ be a $\gamma ^{c}$-function of $T_{p}^{s}$.
Then for $l~(1\leq l\leq s)$, $\sum _{u\in V(L_{p}^{l})}f(u)\geq \min\{cp+1,p\}$.
\end{lem}
%%%%%%%%%%%%%%%%%%%%%%%%%%%%%%%%%%%%%%%%%%%%%%%%%%%%%%%%%%%%%%%%%%%%%%%%%%%%%%%%%%%%%%%%%%%%%%%%%%%%%%%%%%%%%%%%%%%%%%%%
\proof
If $f(x^{l})\geq 1$, then $f(y_{i,2}^{l})\geq c$ or $f(y_{i,1}^{l})\geq 1$ for each $i~(1\leq i\leq p)$, and hence
$$
\sum _{u\in V(L_{p}^{l})}f(u)=f(x^{l})+\sum _{1\leq i\leq p}(f(y_{i,1}^{l})+f(y_{i,2}^{l}))\geq 1+\sum _{1\leq i\leq p}c=1+cp,
$$
as desired.
Thus we may assume that $f(x^{l})<1$.
For each $i~(1\leq i\leq p)$, since the restriction of $f$ on $\{y_{i,1},y_{i,2}\}$ is a $c$-SDF of $G[\{y_{i,1},y_{i,2}\}]~(\simeq P_{2})$, we have $f(y_{i,1})+f(y_{i,2})\geq \gamma ^{c}(P_{2})=1$.
It follows that
$$
\sum _{u\in V(L_{p}^{l})}f(u)=f(x^{l})+\sum _{1\leq i\leq p}(f(y_{i,1}^{l})+f(y_{i,2}^{l}))\geq \sum _{1\leq i\leq p}1=p,
$$
as desired.
\qed

Now we show that Theorem~\ref{mainthm} for the case where $\frac{1}{2}\leq c<1$ is best possible.
We assume that $\frac{m}{m+1}\leq c<\frac{2m+1}{2m+3}$.
Let $f$ be a $\gamma ^{c}$-function of $T_{m+1}^{s}$.
%Fix an index $l~(1\leq l\leq s)$.
Since $c(m+1)+1\geq \frac{m(m+1)}{m+1}+1=m+1$ and $|V(T_{m+1}^{s})|=s(2m+3)$, it follows from Lemma~\ref{lem6.1-1} that
\begin{align*}
\gamma ^{c}(T_{m+1}^{s}) &= w(f)\\
&= \sum _{1\leq l\leq s}\left(\sum _{u\in V(L_{m+1}^{l})}f(u)\right)\\
&\geq \sum _{1\leq l\leq s}\min\{c(m+1)+1,m+1\}\\
&= s(m+1)\\
&= \frac{m+1}{2m+3}|V(T_{m+1}^{s})|.
\end{align*}
This together with Theorem~\ref{mainthm} implies that $\gamma ^{c}(T_{m+1}^{s})=\frac{m+1}{2m+3}|V(T_{m+1}^{s})|$.
Since $s\geq 1$ is arbitrary, there exist infinitely many connected graphs $G$ with $\gamma ^{c}(G)=\frac{m+1}{2m+3}|V(G)|$.

Next we assume that $\frac{2m+1}{2m+3}\leq c<\frac{m+1}{m+2}$.
Let $f'$ be a $\gamma ^{c}$-function of $T_{m+2}^{s}$.
%Fix an index $l~(1\leq l\leq s)$.
Since $c(m+2)+1<\frac{(m+1)(m+2)}{m+2}+1=m+2$ and $|V(T_{m+2}^{s})|=s(2m+5)$, it follows from Lemma~\ref{lem6.1-1} that
\begin{align*}
\gamma ^{c}(T_{m+2}^{s}) &= w(f')\\
&= \sum _{1\leq l\leq s}\left(\sum _{u\in V(L_{m+2}^{l})}f'(u)\right)\\
&\geq \sum _{1\leq l\leq s}\min\{c(m+2)+1,m+2\}\\
&= s(c(m+2)+1)\\
&= \frac{cm+2c+1}{2m+5}|V(T_{m+2}^{s})|.
\end{align*}
This together with Theorem~\ref{mainthm} implies that $\gamma ^{c}(T_{m+2}^{s})=\frac{cm+2c+1}{2m+5}|V(T_{m+2}^{s})|$.
Since $s\geq 1$ is arbitrary, there exist infinitely many connected graphs $G$ with $\gamma ^{c}(G)=\frac{cm+2c+1}{2m+5}|V(G)|$.

Therefore, Theorem~\ref{mainthm} for the case where $\frac{1}{2}\leq c<1$ is best possible.

%%%%%%%%%%%%%%%%%%%%%%%%%%%%%%%%%%%%%%%%%%%%%%%%%%%%%%%%%%%%%%%%%%%%%%%%%%%%%%%%%%%%%%%%%%%%%%%%%%%%%%%%%%%%%%%%%%%%%%%%
%%%%%%%%%%%%%%%%%%%%%%%%%%%%%%%%%%%%%%%%%%%%%%%%%%%%%%%%%%%%%%%%%%%%%%%%%%%%%%%%%%%%%%%%%%%%%%%%%%%%%%%%%%%%%%%%%%%%%%%%
\subsection{The case $c=1$}\label{sec6.2}
%%%%%%%%%%%%%%%%%%%%%%%%%%%%%%%%%%%%%%%%%%%%%%%%%%%%%%%%%%%%%%%%%%%%%%%%%%%%%%%%%%%%%%%%%%%%%%%%%%%%%%%%%%%%%%%%%%%%%%%%
%%%%%%%%%%%%%%%%%%%%%%%%%%%%%%%%%%%%%%%%%%%%%%%%%%%%%%%%%%%%%%%%%%%%%%%%%%%%%%%%%%%%%%%%%%%%%%%%%%%%%%%%%%%%%%%%%%%%%%%%

Fink et al.~\cite{FJKR} and Payan and Xuong~\cite{PX} proved that a connected graph $G$ satisfies $\gamma (G)=\frac{1}{2}|V(G)|$ if and only if $G$ is isomorphic to either a cycle of order $4$ or the graph obtained from a connected graph $H$ by adding a pendant edge to each vertex of $H$.
This together with Proposition~\ref{prop1.1}(i) implies that there exist infinitely many connected graphs $G$ with $\gamma ^{1}(G)=\frac{1}{2}|V(G)|$.
Consequently Theorem~\ref{mainthm} for the case where $c=1$ is best possible.

%%%%%%%%%%%%%%%%%%%%%%%%%%%%%%%%%%%%%%%%%%%%%%%%%%%%%%%%%%%%%%%%%%%%%%%%%%%%%%%%%%%%%%%%%%%%%%%%%%%%%%%%%%%%%%%%%%%%%%%%
%%%%%%%%%%%%%%%%%%%%%%%%%%%%%%%%%%%%%%%%%%%%%%%%%%%%%%%%%%%%%%%%%%%%%%%%%%%%%%%%%%%%%%%%%%%%%%%%%%%%%%%%%%%%%%%%%%%%%%%%
\subsection{The case $1<c\leq 2$}\label{sec6.3}
%%%%%%%%%%%%%%%%%%%%%%%%%%%%%%%%%%%%%%%%%%%%%%%%%%%%%%%%%%%%%%%%%%%%%%%%%%%%%%%%%%%%%%%%%%%%%%%%%%%%%%%%%%%%%%%%%%%%%%%%
%%%%%%%%%%%%%%%%%%%%%%%%%%%%%%%%%%%%%%%%%%%%%%%%%%%%%%%%%%%%%%%%%%%%%%%%%%%%%%%%%%%%%%%%%%%%%%%%%%%%%%%%%%%%%%%%%%%%%%%%

Throughout this subsection, fix a number $c$ with $1<c\leq 2$, and let $m\geq 1$ be the integer such that $\frac{m+2}{m+1}<c\leq \frac{m+1}{m}$.

\begin{lem}%%%%%%%%%%%%%%%%%%%%%%%%%%%%%%%%%%%%%%%%%%%%%%%%%%%%%%%%%%%%%%%%%%%%%%%%%%%%%%%%%%%%%%%%%%%%%%%%%%%%%%%%%%%%%
\label{lem6.3-1}
Let $p\geq 1$ and $s\geq 1$ be integers.
Let $f:V(T_{p}^{s})\rightarrow \mathbb{R}^{+}$ be a $\gamma ^{c}$-function of $T_{p}^{s}$.
Then for $l~(1\leq l\leq s)$, $\sum _{u\in V(L_{p}^{l})}f(u)\geq \min\{p+1,cp\}$.
\end{lem}
%%%%%%%%%%%%%%%%%%%%%%%%%%%%%%%%%%%%%%%%%%%%%%%%%%%%%%%%%%%%%%%%%%%%%%%%%%%%%%%%%%%%%%%%%%%%%%%%%%%%%%%%%%%%%%%%%%%%%%%%
\proof
If $f(x^{l})\geq 1$, then $f(y_{i,2}^{l})\geq c$ or $f(y_{i,1}^{l})\geq 1$ for each $i~(1\leq i\leq p)$, and hence
$$
\sum _{u\in V(L_{p}^{l})}f(u)=f(x^{l})+\sum _{1\leq i\leq p}(f(y_{i,1}^{l})+f(y_{i,2}^{l}))\geq 1+\sum _{1\leq i\leq p}1=1+p,
$$
as desired.
Thus we may assume that $f(x^{l})<1$.
For each $i~(1\leq i\leq p)$, since the restriction of $f$ on $\{y_{i,1},y_{i,2}\}$ is a $c$-SDF of $G[\{y_{i,1},y_{i,2}\}]~(\simeq P_{2})$, we have $f(y_{i,1})+f(y_{i,2})\geq \gamma ^{c}(P_{2})=c$.
It follows that
$$
\sum _{u\in V(L_{p}^{l})}f(u)=f(x^{l})+\sum _{1\leq i\leq p}(f(y_{i,1}^{l})+f(y_{i,2}^{l}))\geq \sum _{1\leq i\leq p}c=cp,
$$
as desired.
\qed

Now we show that Theorem~\ref{mainthm} for the case where $1<c\leq 2$ is best possible.
We assume that $\frac{m+2}{m+1}<c\leq \frac{(2m+1)(m+2)}{(2m+3)m}$.
Let $f$ be a $\gamma ^{c}$-function of $T_{m+1}^{s}$.
%Fix an index $l~(1\leq l\leq s)$.
Since $c(m+1)\geq \frac{(m+2)(m+1)}{m+1}=m+2$ and $|V(T_{m+1}^{s})|=s(2m+3)$, it follows from Lemma~\ref{lem6.3-1} that
\begin{align*}
\gamma ^{c}(T_{m+1}^{s}) &= w(f)\\
&= \sum _{1\leq l\leq s}\left(\sum _{u\in V(L_{m+1}^{l})}f(u)\right)\\
&\geq \sum _{1\leq l\leq s}\min\{(m+1)+1,c(m+1)\}\\
&= s(m+2)\\
&= \frac{m+2}{2m+3}|V(T_{m+1}^{s})|.
\end{align*}
This together with Theorem~\ref{mainthm} implies that $\gamma ^{c}(T_{m+1}^{s})=\frac{m+2}{2m+3}|V(T_{m+1}^{s})|$.
Since $s\geq 1$ is arbitrary, there exist infinitely many connected graphs $G$ with $\gamma ^{c}(G)=\frac{m+2}{2m+3}|V(G)|$.

Next we assume that $\frac{(2m+1)(m+2)}{(2m+3)m}<c\leq \frac{m+1}{m}$.
Let $f'$ be a $\gamma ^{c}$-function of $T_{m}^{s}$.
%Fix an index $l~(1\leq l\leq s)$.
Since $cm\leq \frac{(m+1)m}{m}=m+1$ and $|V(T_{m}^{s})|=s(2m+1)$, it follows from Lemma~\ref{lem6.3-1} that
\begin{align*}
\gamma ^{c}(T_{m}^{s}) &= w(f')\\
&= \sum _{1\leq l\leq s}\left(\sum _{u\in V(L_{m}^{l})}f'(u)\right)\\
&\geq \sum _{1\leq l\leq s}\min\{m+1,cm\}\\
&= scm\\
&= \frac{cm}{2m+1}|V(T_{m}^{s})|.
\end{align*}
This together with Theorem~\ref{mainthm} implies that $\gamma ^{c}(T_{m}^{s})=\frac{cm}{2m+1}|V(T_{m}^{s})|$.
Since $s\geq 1$ is arbitrary, there exist infinitely many connected graphs $G$ with $\gamma ^{c}(G)=\frac{cm}{2m+1}|V(G)|$.

Therefore, Theorem~\ref{mainthm} for the case where $1<c\leq 2$ is best possible.

%%%%%%%%%%%%%%%%%%%%%%%%%%%%%%%%%%%%%%%%%%%%%%%%%%%%%%%%%%%%%%%%%%%%%%%%%%%%%%%%%%%%%%%%%%%%%%%%%%%%%%%%%%%%%%%%%%%%%%%%
%%%%%%%%%%%%%%%%%%%%%%%%%%%%%%%%%%%%%%%%%%%%%%%%%%%%%%%%%%%%%%%%%%%%%%%%%%%%%%%%%%%%%%%%%%%%%%%%%%%%%%%%%%%%%%%%%%%%%%%%
\subsection{The case $c>2$}\label{sec6.4}
%%%%%%%%%%%%%%%%%%%%%%%%%%%%%%%%%%%%%%%%%%%%%%%%%%%%%%%%%%%%%%%%%%%%%%%%%%%%%%%%%%%%%%%%%%%%%%%%%%%%%%%%%%%%%%%%%%%%%%%%
%%%%%%%%%%%%%%%%%%%%%%%%%%%%%%%%%%%%%%%%%%%%%%%%%%%%%%%%%%%%%%%%%%%%%%%%%%%%%%%%%%%%%%%%%%%%%%%%%%%%%%%%%%%%%%%%%%%%%%%%

In Subsection~\ref{sec6.3}, we show that there exist infinitely many connected graphs $G$ with $\gamma ^{2}(G)=\frac{2}{3}|V(G)|$.
On the other hand, it follows from Lemma~\ref{lem5.1} that for a number $c>2$ and a connected graph $G$ of order at least $2$, $\gamma ^{c}(G)=\gamma ^{2}(G)$.
Hence for $c>2$, there exist infinitely many connected graphs $G$ with $\gamma ^{c}(G)=\frac{2}{3}|V(G)|$.
Consequently Theorem~\ref{mainthm} for the case where $c>2$ is best possible.

%%%%%%%%%%%%%%%%%%%%%%%%%%%%%%%%%%%%%%%%%%%%%%%%%%%%%%%%%%%%%%%%%%%%%%%%%%%%%%%%%%%%%%%%%%%%%%%%%%%%%%%%%%%%%%%%%%%%%%%%
%%%%%%%%%%%%%%%%%%%%%%%%%%%%%%%%%%%%%%%%%%%%%%%%%%%%%%%%%%%%%%%%%%%%%%%%%%%%%%%%%%%%%%%%%%%%%%%%%%%%%%%%%%%%%%%%%%%%%%%%
%%%%%%%%%%%%%%%%%%%%%%%%%%%%%%%%%%%%%%%%%%%%%%%%%%%%%%%%%%%%%%%%%%%%%%%%%%%%%%%%%%%%%%%%%%%%%%%%%%%%%%%%%%%%%%%%%%%%%%%%
\section{Concluding remarks}\label{sec-cm}
%%%%%%%%%%%%%%%%%%%%%%%%%%%%%%%%%%%%%%%%%%%%%%%%%%%%%%%%%%%%%%%%%%%%%%%%%%%%%%%%%%%%%%%%%%%%%%%%%%%%%%%%%%%%%%%%%%%%%%%%
%%%%%%%%%%%%%%%%%%%%%%%%%%%%%%%%%%%%%%%%%%%%%%%%%%%%%%%%%%%%%%%%%%%%%%%%%%%%%%%%%%%%%%%%%%%%%%%%%%%%%%%%%%%%%%%%%%%%%%%%
%%%%%%%%%%%%%%%%%%%%%%%%%%%%%%%%%%%%%%%%%%%%%%%%%%%%%%%%%%%%%%%%%%%%%%%%%%%%%%%%%%%%%%%%%%%%%%%%%%%%%%%%%%%%%%%%%%%%%%%%

In this paper, our aim is to find essential boundaries among domination, total domination and Roman domination.
Thus we focused on $c$-self-domination for $c\geq \frac{1}{2}$.
As far as we check Theorem~\ref{mainthm} and its proof, the essential parts are roughly divided into $\frac{1}{2}\leq c<1$, $c=1$ and $c>1$.
Hence, for example, we expect that some results on Roman domination can be extended to $c$-self-domination for $\frac{1}{2}<c<1$.

Here one might be interested in upper bounds on $c$-self-domination for $c<\frac{1}{2}$.
As a general upper bound on $c$-self-domination, we obtain the following proposition.

\begin{prop}%%%%%%%%%%%%%%%%%%%%%%%%%%%%%%%%%%%%%%%%%%%%%%%%%%%%%%%%%%%%%%%%%%%%%%%%%%%%%%%%%%%%%%%%%%%%%%%%%%%%%%%%%%%%
\label{prop-remark1}
Let $0<c\leq 1$ be a number, and let $G$ be a graph of order $n$.
Then $\gamma ^{c}(G)\leq cn$.
\end{prop}
%%%%%%%%%%%%%%%%%%%%%%%%%%%%%%%%%%%%%%%%%%%%%%%%%%%%%%%%%%%%%%%%%%%%%%%%%%%%%%%%%%%%%%%%%%%%%%%%%%%%%%%%%%%%%%%%%%%%%%%%
\proof
Since the function $f:V(G)\rightarrow \mathbb{R}^{+}$ with $f(u)=c~(u\in V(G))$ is a $c$-SDF of $G$ with $w(f)=cn$, we get the desired conclusion.
\qed

Proposition~\ref{prop-remark1} is best possible for the case where $0<c\leq \frac{1}{3}$.
It suffices to show that $\gamma ^{c}(P_{n})\geq cn$ for all $n\geq 1$.
Let $f$ be a $\gamma ^{c}$-function of $P_{n}$.
By (\ref{01c}), we have $\{f(u):u\in V(P_{n})\}\subseteq \{0,1,c\}$.
Choose $f$ so that $|\{u\in V(P_{n}):f(u)=c\}|$ is as large as possible.
Suppose that $\{f(u):u\in V(P_{n})\}\cap \{0,1\}\neq \emptyset $.
Then there exists a vertex $x\in V(P_{n})$ with $f(x)=1$.
Let $f':V(P_{n})\rightarrow \mathbb{R}^{+}$ be the function with
$$
f'(u)=
\begin{cases}
c, & u=x;\\
\max\{f(u),c\}, & u\in N_{P_{n}}(x);\\
f(u), & \mbox{otherwise}.
\end{cases}
$$
Then $f'$ is a $c$-SDF of $P_{n}$ with $w(f')\leq w(f)-1+3c\leq w(f)$ and $|\{u\in V(P_{n}):f'(u)=c\}|>|\{u\in V(P_{n}):f(u)=c\}|$, which contradicts the choice of $f$.
Thus $f(u)=c$ for all $u\in V(P_{n})$.
Consequently, $\gamma ^{c}(P_{n})=cn$, and so Proposition~\ref{prop-remark1} is best possible for the case where $0<c\leq \frac{1}{3}$.

The author recently settled the remaining case in the following paper~\cite{F}, i.e., for a number $c\in \mathbb{R}^{+}$ with $\frac{1}{3}<c<\frac{1}{2}$, a sharp upper bound on $\gamma ^{c}(G)$ for a connected graph $G$ is given.

%%%%%%%%%%%%%%%%%%%%%%%%%%%%%%%%%%%%%%%%%%%%%%%%%%%%%%%%%%%%%%%%%%%%%%%%%%%%%%%%%%%%%%%%%%%%%%%%%%%%%%%%%%%%%%%%%%%%%%%%
%%%%%%%%%%%%%%%%%%%%%%%%%%%%%%%%%%%%%%%%%%%%%%%%%%%%%%%%%%%%%%%%%%%%%%%%%%%%%%%%%%%%%%%%%%%%%%%%%%%%%%%%%%%%%%%%%%%%%%%%
%%%%%%%%%%%%%%%%%%%%%%%%%%%%%%%%%%%%%%%%%%%%%%%%%%%%%%%%%%%%%%%%%%%%%%%%%%%%%%%%%%%%%%%%%%%%%%%%%%%%%%%%%%%%%%%%%%%%%%%%
\section*{Acknowledgment}
%%%%%%%%%%%%%%%%%%%%%%%%%%%%%%%%%%%%%%%%%%%%%%%%%%%%%%%%%%%%%%%%%%%%%%%%%%%%%%%%%%%%%%%%%%%%%%%%%%%%%%%%%%%%%%%%%%%%%%%%
%%%%%%%%%%%%%%%%%%%%%%%%%%%%%%%%%%%%%%%%%%%%%%%%%%%%%%%%%%%%%%%%%%%%%%%%%%%%%%%%%%%%%%%%%%%%%%%%%%%%%%%%%%%%%%%%%%%%%%%%
%%%%%%%%%%%%%%%%%%%%%%%%%%%%%%%%%%%%%%%%%%%%%%%%%%%%%%%%%%%%%%%%%%%%%%%%%%%%%%%%%%%%%%%%%%%%%%%%%%%%%%%%%%%%%%%%%%%%%%%%

The author would like to thank anonymous referees for careful reading and helpful comments.
This work was supported by JSPS KAKENHI Grant number 18K13449.

\end{document}